\documentclass[12pt,reqno,a4paper]{amsart}
\usepackage{fullpage}
\usepackage{amsmath,amssymb,amsthm,mathtools}
\usepackage{booktabs}
\usepackage[T1]{fontenc}
\usepackage{hyperref}
\hypersetup{hidelinks}
\newcommand{\Q}{\mathbb Q}
\newcommand{\Z}{\mathbb Z}
\newcommand{\C}{\mathbb C}
\newcommand{\R}{\mathbb R}
\newcommand{\B}{\mathbb B}
\newcommand{\PP}{\mathbb P}
\newcommand{\Ad}{\operatorname{Ad}}
\newcommand{\tr}{\operatorname{tr}}
\newcommand{\Tr}{\operatorname{Tr}}
\newcommand{\PU}{\operatorname{PU}}
\newcommand{\Stab}{\operatorname{Stab}}
\newcommand{\Fix}{\operatorname{Fix}}
\newcommand{\Comm}{\operatorname{Comm}}
\newcommand{\chiorb}{\chi_{\mathrm{orb}}}
\newcommand{\kAd}{k_{\mathrm{Ad}}}
\newcommand{\ov}[1]{\overline{#1}}
\newtheorem{thm}{Theorem}[section]
\newtheorem{lem}[thm]{Lemma}
\newtheorem{prop}[thm]{Proposition}

\theoremstyle{remark}
\newtheorem{rmk}[thm]{Remark}
\title[A nonarithmetic lattice from a mirror]{A nonarithmetic lattice from a mirror in the Couwenberg--Heckman--Looijenga construction}
\author{Chenglong Yu}
\date{October 29, 2026}
\newcommand{\Addresses}{{%
 \bigskip
 \footnotesize
 C.~Yu, \textsc{Center for Mathematics and Interdisciplinary Sciences, Fudan University and Shanghai Institute for Mathematics and Interdisciplinary Sciences (SIMIS), Shanghai, China}\par\nopagebreak
 \textit{Email address}: \texttt{yuchenglong@simis.cn}
}}
\begin{document}
\bibliographystyle{amsalpha}
\begin{abstract}
We prove that the stabilizer of a reflection mirror of the
Couwenberg--Heckman--Looijenga lattice $\mathcal C(G_{29},3)$ induces
a nonarithmetic lattice in $\PU(2,1)$. This lattice is not commensurable with any
Deligne--Mostow lattice in $\PU(2,1)$ or any lattice in the
twenty-two nonarithmetic classes listed by Deraux, Parker and Paupert.
\end{abstract}
\maketitle

\section{Introduction}

Margulis' superrigidity theorem \cite{Margulis91} implies that every
lattice in a connected simple real algebraic group of real rank at
least two is arithmetic \cite{Margulis84}. The same conclusion holds for
$\operatorname{Sp}(n,1)$, $n\ge2$, and $F_{4(-20)}$ by the
superrigidity theorems of Corlette \cite{Corlette92} and
Gromov--Schoen \cite[Theorem~8.4]{GS92}. Thus, among connected simple
real algebraic groups with noncompact real points, nonarithmetic
lattices can occur only in the real and complex hyperbolic cases.

In real hyperbolic space, Gromov and Piatetski-Shapiro \cite{GPS87}
constructed infinitely many commensurability classes of nonarithmetic
lattices in every dimension at least two by gluing pieces of arithmetic
manifolds along totally geodesic hypersurfaces. For $\PU(n,1)$,
whether there are infinitely many such classes remains open for every
fixed $n\ge2$. Known examples occur only in complex dimensions two and
three; see \cite[Questions~3.4--3.5]{Fisher23}.

Mostow \cite[Theorem~A]{Mostow80} constructed the first nonarithmetic
lattices in $\PU(2,1)$ using complex reflections and fundamental
polyhedra. Building on Picard's work \cite{Picard81,Picard85},
Deligne and Mostow \cite{DM86} established lattice criteria for
hypergeometric monodromy groups, which Mostow
\cite{Mostow86,Mostow88} subsequently extended.
The resulting ball quotients arise from moduli spaces of
weighted points on $\PP^1$. Thurston \cite[Theorem~0.2]{Thurston98}
gave a geometric description through Euclidean cone metrics on the
sphere. Commensurability relations among the resulting lattices were
established by Sauter \cite{Sauter90} and Deligne--Mostow \cite{DM93}.
Together with the invariants of Kappes--M\"oller \cite{KM16} and
McMullen \cite{McM17}, these results give nine nonarithmetic
Deligne--Mostow commensurability classes in $\PU(2,1)$; see also
\cite[Table~2]{YuZheng}.

Deraux, Parker and Paupert \cite{DPP16} constructed five further
classes. Their subsequent work \cite[Theorem~1.4 and Appendix~A.11]{DPP}
adds eight classes and gives a catalogue of twenty-two: the nine
Deligne--Mostow classes and thirteen classes represented by sporadic
and Thompson triangle groups.

Another approach uses the Miyaoka--Yau inequality.
Yau \cite[Theorem~4]{Yau77} proved that a smooth projective surface
with ample canonical bundle and $c_1^2=3c_2$ has universal cover $\B^2$.
Hirzebruch \cite[\S\S3.1--3.2]{Hirz83} used this criterion to construct
compact ball quotients by resolving abelian covers of $\PP^2$ branched over line
arrangements, including the complete quadrangle and the Hesse
arrangement. Barthel, Hirzebruch and H\"ofer \cite[Chs.~4--5]{BHH87} developed
this method systematically, describing ball quotients through line
configurations and their ramification data.

Couwenberg, Heckman and Looijenga \cite{CHL} gave a common extension
of the Deligne--Mostow and line-arrangement constructions using complex
reflection arrangements. In dimension three, the Deligne--Mostow list
contains one nonarithmetic commensurability class. Deraux
\cite[Theorem~1.2 and \S6.3]{DerauxLattice} proved that
$\Gamma_{29}=\mathcal C(G_{29},3)\subset\PU(3,1)$ gives a second class:
it is nonarithmetic and is not commensurable with a Deligne--Mostow
lattice. Its quotient $\Gamma_{29}\backslash\B^3$ is not compact.

Longitudinal restriction is part of the Dunkl-system theory of
\cite[Lemma~2.18, Definition~2.20 and Lemma~4.4]{CHL}.
The restricted arrangement $(G_{29},A_1)$ has also been studied
\cite[Table~3]{AHR}, \cite[\S1]{ADR}.
Deraux \cite[\S\S1,10]{DerauxMirrors} studied the Fuchsian stabilizers of
reflection mirrors in lattice complex hyperbolic triangle groups
in $\PU(2,1)$.

A two-ball group associated with an exceptional divisor already
appears in Deraux's initial preprint \cite[Figure~7]{DerauxEarly}.
The rank-three transversal system at a point where twelve mirrors
meet gives Mostow's $\Gamma(3,1/3)$.

Let $M\cong\B^2$ be the mirror of an order-three reflection in
$\Gamma_{29}$ corresponding, under the Couwenberg--Heckman--Looijenga
construction, to a reflecting hyperplane of $G_{29}$.
All such mirrors are conjugate. Define
\[
 \Lambda_{29}=\Stab_{\Gamma_{29}}(M)/\Fix_{\Gamma_{29}}(M)
 \subset\PU(2,1),
\]
where $\Fix_{\Gamma_{29}}(M)$ is the pointwise stabilizer.

For a lattice $\Lambda\subset\PU(n,1)$, write
$\chiorb(\Lambda)=\chiorb(\Lambda\backslash\B^n)$ for its orbifold
Euler characteristic and
$\kAd(\Lambda)=\Q(\Tr\Ad(g):g\in\Lambda)$ for its adjoint trace field.

\begin{thm}\label{thm:main}
The group $\Lambda_{29}$ is a nonarithmetic lattice, and
$\Lambda_{29}\backslash\B^2$ is not compact. The
adjoint trace field and orbifold Euler characteristic are
\[
 \kAd(\Lambda_{29})=\Q(\sqrt3),\qquad
 \chiorb(\Lambda_{29}\backslash\B^2)=\frac{43}{72}.
\]
Every elliptic boundary curve of a neat finite cover is isogenous to
$\C/\Z[\omega]$, where $\omega=e^{2\pi i/3}$.
Moreover, $\Lambda_{29}$ is not commensurable with any Deligne--Mostow
lattice in $\PU(2,1)$ or with any lattice in the twenty-two
nonarithmetic classes listed in \cite[Appendix~A.11]{DPP}.
\end{thm}

The quotient $\Lambda_{29}\backslash M$ is the normalization of the
image of $M$ in $\Gamma_{29}\backslash\B^3$, with its induced orbifold
structure. In the finite arrangement model,
the map from $S=\operatorname{Bl}_6\PP^2$ to the compactified mirror
quotient has degree $24$.
The orbifold canonical class pulls back to
$7h-\sum_{j=1}^6\varepsilon_j$, of square $43$, where $h$ is the
plane class and the $\varepsilon_j$ are the exceptional classes.
Extending the canonical curvature identity across the cusps and applying
the projection formula gives the Euler characteristic $43/72$.
We also compute it by stratification.

The trace field and the fact that $\Lambda_{29}\backslash\B^2$ is not
compact leave two possible classes in the catalogue of Deraux, Parker
and Paupert. The first is
excluded by comparing the isogeny classes of boundary elliptic curves,
using Yu--Zheng \cite[Proposition~8.1]{YuZheng}.
The second is represented by Mostow's
$\Gamma(6,1/6)$, with Euler characteristic $11/144$.
For this class we follow the proof of \cite[Proposition~7.9]{DPP}.
Multiplicativity of the orbifold Euler characteristic would force the
index of this Mostow group in a common lattice overgroup to be divisible
by $11$. The cusp estimate in \cite[\S7.2.2, item~(3) following
Proposition~7.16]{DPP} bounds that index by three.

Section~\ref{sec:finite-volume} identifies the mirror quotient.
Section~\ref{sec:matrices} computes its trace field and cusp type and
proves nonarithmeticity. Sections~\ref{sec:volume} and~\ref{sec:euler}
compute the orbifold Euler characteristic, and
Section~\ref{sec:comparison} compares $\Lambda_{29}$ with the classes
in Table~\ref{tab:dpp-classes} and proves Theorem~\ref{thm:comparison}.

All varieties are over $\C$. We use negative lines for Hermitian forms of signature $(n,1)$.
Matrices below are written in unitary lifts; their projective
classes define the groups.

\section{The mirror quotient}
\label{sec:finite-volume}

We recall the fixed-point argument for the lattice property of a
mirror stabilizer.

\begin{lem}\label{lem:mirror}
Let \(\Gamma\subset\PU(n,1)\), \(n\ge2\), be a lattice, and let
\(r\in\Gamma\) be a nontrivial finite-order complex reflection with
mirror \(M\cong\B^{n-1}\). The image of \(\Stab_\Gamma(M)\)
in \(\PU(n-1,1)\) is a lattice.
\end{lem}

\begin{proof}
Since $\Gamma$ is finitely generated and linear, Selberg's lemma
gives a torsion-free subgroup of finite index. Taking its normal
core gives a torsion-free normal finite-index subgroup $\Gamma_0$.
Let \(f\) be the finite-order isometry induced by \(r\)
on \(X_0=\Gamma_0\backslash\B^n\). The image of \(M\) is a connected
component \(C\) of \(\Fix(f)\).

If \(g\in\Gamma_0\)
and \(gM\cap M\ne\varnothing\), then
\[
 grg^{-1}r^{-1}\in\Gamma_0
\]
fixes an interior point. A discrete torsion-free group has no nontrivial
elliptic element, so this commutator is the identity and \(gM=M\).
More generally, two lifts of \(f\) fixing a common interior point agree.
The local fixed-locus description and continuation along paths therefore
identify
\[
 C=M/(\Gamma_0\cap\Stab_\Gamma(M))
\]
with an embedded connected fixed component. It is closed.

Choose a cusp truncation of \(X_0\) invariant under the finite group
\(\Gamma/\Gamma_0\). The intersection of \(C\) with the truncated
compact core is compact. If a cusp contains a point of \(C\), it is
preserved by \(f\). A lift preserving its ideal point has dilation one:
a finite power lies in the cusp group, which has trivial dilation.
Consequently it preserves the height coordinate and acts on the
horosphere quotient independently of height.
Its fixed locus in that compact horosphere quotient is closed and
compact. There are finitely many such components, and the height flow
preserves them. Each resulting end of \(C\) is an induced
complex hyperbolic cusp with compact cross-section. In complex dimension
\(m=n-1\ge1\), its volume density decays exponentially in height
(as \(e^{-2mt}\) in a standard normalization), so its volume is finite.
Thus \(C\) has finite volume.

The restriction map from \(\Stab_{\PU(n,1)}(M)\) to \(\PU(n-1,1)\)
has compact kernel, so the image of \(\Stab_\Gamma(M)\) is discrete.
Its quotient is a finite quotient of \(C\), hence has finite volume.
\end{proof}

\subsection{Normalization of the mirror divisor}

We identify \(\Lambda_{29}\backslash M\) with the normalization of the
image of \(M\) in \(\Gamma_{29}\backslash\B^3\).

We use the standard normalization property for finite quotients.
Let a finite group \(F\) act on \(X_0\), and let \(C\subset X_0\) be
a normal irreducible closed subvariety. The normalization of its image
in \(X_0/F\) is
\[
 C/\Stab_F(C).
\]
Indeed, the quotient is normal and its map to the image of \(C\) is
finite and birational. The same argument applies to closed analytic
subvarieties. For $\Gamma=\Gamma_{29}$, apply this with
$F=\Gamma/\Gamma_0$ to the embedded fixed component above.
The preimage of \(\Stab_F(C)\) in \(\Gamma\) is
\(\Gamma_0\Stab_\Gamma(M)\), which identifies the normalization with
the mirror quotient.

Write \(PG_{29}=G_{29}/Z(G_{29})\). Apply the same argument to
the noncompact Couwenberg--Heckman--Looijenga model \(Y^0/PG_{29}\) of
\cite[Proposition~2(2)]{DerauxVolumes}, where \(Y^0\) is obtained from the
compact blowup model by deleting the cusp points. The preimage of this
branch divisor is the union of the strict transforms of the forty
projectivized reflecting hyperplanes of \(G_{29}\).
If \(S\) is one compact strict transform and
\(S^0=S\cap Y^0\), the normalization of the open branch divisor has
coarse space \(S^0/\Stab_{PG_{29}}(S)\).
We equip this space with the orbifold structure induced by
$\Lambda_{29}$ in the $p=3$ ball quotient. The finite $p=2$
stabilizer orders enter the calculations through their ratios with
the stabilizer orders for this deformed action.

\section{The trace field and the cusp}
\label{sec:matrices}

Put
\[
 \omega=e^{2\pi i/3},\qquad
 a=\frac1{\omega-1}=-\frac{2+\omega}{3},\qquad
 b=\frac{1+i}{\omega(\omega-1)}
   =\frac{(1+i)(1+2\omega)}3.
\]
Deraux's realization \cite[\S5.2, (5-8)--(5-9), Proposition~5.4]{DerauxLattice}
has Hermitian form
\[
 H_4=\begin{pmatrix}
 1&a&0&0\\
 \ov a&1&a&b\\
 0&\ov a&1&a\\
 0&\ov b&\ov a&1
 \end{pmatrix}
\]
and generators $R_j=I+(\omega-1)e_je_j^*H_4$, $1\le j\le4$,
where $e_1,\ldots,e_4$ are the standard basis vectors.
Each generator is the identity except for its $j$-th row, which is
\begin{equation}\label{eq:ambient-rows}
\begin{aligned}
 R_1&:\ (\omega,1,0,0),&
 R_2&:\ (-\omega,\omega,1,(1+i)\omega^2),\\
 R_3&:\ (0,-\omega,\omega,1),&
 R_4&:\ (0,-(1-i)\omega^2,-\omega,\omega).
\end{aligned}
\end{equation}
Thus $\Gamma_{29}=\langle[R_1],[R_2],[R_3],[R_4]\rangle$.
The parameter here is $\mu=1+i$, the discrete branch in
\cite[Proposition~5.4]{DerauxLattice}.
We take $M=\PP(e_1^\perp)\cap\B^3$.

In the ordered basis
\[
 e_2-ae_1,\quad e_3,\quad e_4
\]
of \(e_1^\perp\), the induced Hermitian form is
\begin{equation}\label{eq:mirror-form}
 G=
 \begin{pmatrix}
 2/3&a&b\\
 \ov a&1&a\\
 \ov b&\ov a&1
 \end{pmatrix}.
\end{equation}
Its leading principal minors are
\[
 \frac23,\quad\frac13,\quad-\frac{2+\sqrt3}{9},
\]
so its signature is \((2,1)\).
The following ambient elements preserve \(e_1^\perp\) and restrict to
\[
\begin{aligned}
 T=(R_1R_2)^3|_M
 &=\begin{pmatrix}
 -1&2+\omega&-(1+i)(1+2\omega)\\0&1&0\\0&0&1
 \end{pmatrix},\\
 P=R_3|_M
 &=\begin{pmatrix}1&0&0\\-\omega&\omega&1\\0&0&1\end{pmatrix},\\
 Q=R_4|_M
 &=\begin{pmatrix}
 1&0&0\\0&1&0\\-(1-i)\omega^2&-\omega&\omega
 \end{pmatrix}.
\end{aligned}
\]
Direct multiplication gives
\[
 T^2=P^3=Q^3=I,\qquad
 T^*GT=P^*GP=Q^*GQ=G.
\]
We use \(T,P,Q\) as elements of \(\Lambda_{29}\), without assuming
that they generate it.

\subsection{The adjoint trace field}

For a unitary lift of dimension three,
\[
 \Tr\Ad([g])=|\tr g|^2-1.
\]
Put $K=\Q(i,\omega)=\Q(\zeta_{12})$, where $\zeta_{12}=e^{\pi i/6}$.
Every ambient word has entries in $K$. An ambient word preserving
\(e_1^\perp\) restricts in the displayed basis to a matrix over $K$
preserving \(G\). Thus all adjoint traces of \(\Lambda_{29}\) lie in
the real subfield \(\Q(\sqrt3)\). Conversely,
\[
 \tr(TPQ)=3+\sqrt3,\qquad
 \Tr\Ad(TPQ)=11+6\sqrt3.
\]
Consequently
\[
 \kAd(\Lambda_{29})=\Q(\sqrt3).
\]

\subsection{Cusp groups}

Direct multiplication gives
\begin{equation}\label{eq:unipotent}
 (TQ)^6=I+N,\qquad
 N=
 \begin{pmatrix}
 0&-3\sqrt3-i(9+6\sqrt3)&0\\
 0&0&0\\
 0&3+3\sqrt3-i(9+3\sqrt3)&0
 \end{pmatrix}.
\end{equation}
Here \(N\ne0\), \(N^2=0\), and \(\operatorname{rank}N=1\);
its image is isotropic for \(G\). Hence
\(\Lambda_{29}\backslash\B^2\) is not compact.
Moreover
\[
 QN=NQ=N,\qquad TN=N.
\]
If \(v\) spans \(\operatorname{im}N\), then \(Qv=v\).
The third basis vector lies in \(v^\perp\), has positive norm, and is a
\(\omega\)-eigenvector for \(Q\). Therefore the horizontal action at
this cusp, on \(v^\perp/\C v\), is multiplication by \(\omega\).

To describe a Heisenberg lattice at this cusp, put
\[
 c=(1-\omega,0,0)^t,\qquad
 q=(1-\omega,0,-(1-i)\omega^2)^t,
\]
in the basis used in \eqref{eq:mirror-form}. Then $q$ spans $\operatorname{im}N$,
$Tq=Qq=q$, and $q^\perp$ has basis $(q,c)$. Write
\[
 Z=(TQ)^6,\qquad U=[T,Q],\qquad V=QUQ^{-1},
 \qquad [g,h]=ghg^{-1}h^{-1}.
\]
Direct multiplication gives
\begin{gather*}
 [U,V]=Z^{-1},\qquad [Z,T]=[Z,Q]=I,\\
 TUT=U^{-1},\qquad TVT=ZV^{-1},\qquad
 QVQ^{-1}=V^{-1}U^{-1},\\
 Uc=c+2(1-\omega^2)q,\qquad
 Vc=c+2(1-\omega^2)\omega^2q.
\end{gather*}
The matrices $U,V$ are unipotent. Their horizontal parameters are
linearly independent over $\mathbb R$, and $Z$ is a nontrivial central
unipotent. The relation $[U,V]=Z^{-1}$ puts each element of
$\langle U,V\rangle$ in the form $Z^kU^mV^n$.
Horizontal independence and the infinite order of $Z$ show that this
form is unique. It follows that $\langle U,V\rangle$ is a lattice in
the three-dimensional Heisenberg group, with compact quotient.
The displayed conjugation
relations make it normal in $\langle T,Q\rangle$.
The quotient is an abelian quotient of $C_2\times C_3$.
Its horizontal image, generated by $-1$ and $\omega$, has order six,
so the quotient is cyclic of order six.

The cusp stabilizer is discrete and contains
$\langle U,V\rangle$ with finite index. The horizontal translation
lattice of $\langle U,V\rangle$ is a complex multiple of $\Z[\omega]$, so the
elliptic boundary curve in a neat finite cover is isogenous to
\[
 E_\omega=\mathbb C/\mathbb Z[\omega].
\]
In fact every cusp has this isogeny type. The ambient quotient has
one cusp, and its horizontal translation lattice is
$\mathbb Z[\omega]^2$ by \cite[\S6.3, Proposition~6.7]{DerauxLattice}.
Move the ambient cusp containing a given mirror cusp to the standard
one. Intersect the ambient setwise stabilizer of the mirror with a
neat ambient cusp subgroup of finite index, and restrict the intersection
to the mirror. The resulting group has finite index in the mirror cusp
stabilizer, and its translations are restrictions of ambient unipotent
translations. A mirror cusp has
horizontal directions in a complex line $W\subset\mathbb C^2$.
The rational span of its translation lattice lies in
$W\cap\mathbb Q(\omega)^2$. Choose a nonzero vector $v$ in this span.
Then $W=\mathbb C v$ and
\[
 W\cap\mathbb Q(\omega)^2=v\mathbb Q(\omega).
\]
Both the mirror translation lattice and $v\mathbb Z[\omega]$ are full
$\mathbb Z$-lattices in this two-dimensional rational vector space,
so they are commensurable. Thus all toroidal boundary elliptic curves
of neat finite covers are isogenous to $E_\omega$.

\subsection{Nonarithmeticity}

An arithmetic complex hyperbolic lattice whose quotient is not compact
has adjoint trace field \(\Q\). Indeed, the adjoint trace field of an arithmetic lattice
is totally real, and its defining adjoint group is compact at every
nonidentity real place; see Mostow \cite[\S4, p.~192]{Mostow80}.
A nontrivial unipotent has a nontrivial power in an arithmetic
finite-index subgroup. This power remains nontrivial under any field
embedding, whereas a compact group contains no nontrivial unipotents. Thus there
can be no nonidentity real place. The existence of
\eqref{eq:unipotent} and \(\kAd=\Q(\sqrt3)\) therefore prove that
\(\Lambda_{29}\) is nonarithmetic. The conjugate form is also
indefinite. If \(\sigma\) extends the nontrivial embedding of
\(\Q(\sqrt3)\) to \(\Q(i,\omega)\), then
\[
 \det G^\sigma=-\frac{2-\sqrt3}{9}<0.
\]

\section{The mirror surface and its canonical class}
\label{sec:volume}

\subsection{The finite quotient}
\label{sec:finite-quotient}

Specializing Deraux's realization \cite[\S5.2]{DerauxLattice} to
$p=2$ gives generators $r_j$ of $G_{29}$, equal to the identity
except for their $j$-th rows
\[
 (-1,1,0,0),\quad (1,-1,1,-1-i),\quad
 (0,1,-1,1),\quad (0,-1+i,1,-1).
\]
They preserve the Hermitian form
\[
 H_0=\begin{pmatrix}
 1&-1/2&0&0\\
 -1/2&1&-1/2&(1+i)/2\\
 0&-1/2&1&-1/2\\
 0&(1-i)/2&-1/2&1
 \end{pmatrix}.
\]

The orbit of $[e_1]$ consists of the forty mirror polars.
It is computed by applying the $r_j$ and normalizing each vector
by its first nonzero coordinate. The reflecting hyperplane of $r_1$
is parametrized by
\[
 x=(y_1/2,y_1,y_2,y_3).
\]
Restricting $v^*H_0x=0$ to this hyperplane for the other thirty-nine
polars gives twenty-one distinct lines. The number of equations
defining each line, plus one for the chosen mirror, is its ambient
multiplicity. Cross products of line coefficient vectors give all
intersection points. Substitution into the twenty-one line equations
and the forty hyperplane equations gives the line incidences and
ambient multiplicity of each point.

Let \(L_I\) denote the projectivized intersection of the reflecting
hyperplanes of \(r_j\) for \(j\in I\).
For \(p=3\), the Couwenberg--Heckman--Looijenga model of
\cite[\S5.3 and Figure~10]{DerauxVolumes} is
\[
 Y=\operatorname{Bl}_{20}\mathbf P^3 .
\]
The centers form the orbit of \(L_{234}\), and \(Y^0\) is obtained
by deleting the forty points in the \(L_{124}\) cusp orbit.
Each center lies on twelve mirrors. Since the forty mirrors form
one orbit, each mirror contains
\(20\cdot12/40=6\) centers, so its strict transform is
\[
 S=\operatorname{Bl}_{6}\mathbf P^2.
\]
Each cusp point lies on nine mirrors, so \(S^0\) is obtained from
\(S\) by deleting \(40\cdot9/40=9\) cusp points.

The finite group \(G_{29}\) has order \(7680\) and scalar center of order
four \cite[Table~VII]{ST54}. Its projective setwise mirror stabilizer
therefore has order
\[
 \frac{7680}{4\cdot40}=48.
\]
The kernel of its action on the projective mirror has order two.
Indeed, all matrices of the finite \(p=2\) model have entries in
\(\Q(i)\) and determinant in \(\{1,-1\}\).
If \(g\in G_{29}\) acts projectively trivially on a \(\Q(i)\)-rational mirror
hyperplane \(H\), then \(g|_H=\lambda I\) for a root of unity
\(\lambda\in\Q(i)\), hence $\lambda\in\mu_4=\{1,-1,i,-i\}$.
The scalar center is \(\mu_4 I\), so \(\lambda^{-1}g\) fixes \(H\)
linearly. Its remaining eigenvalue is its determinant, namely \(1\)
or \(-1\). It is either the identity or the unique involutory unitary
reflection in \(H\). This proves the assertion.

The quotient map from \(S\) therefore has degree
\[
 \frac{48}{2}=24.
\]
Closing the identity under multiplication by the $r_j$ gives all
$7680$ matrices. Of these, $192$ preserve $\mathbf Ce_1$, and $8$
restrict to scalar matrices on $e_1^{\perp_{H_0}}$. Thus the linear
stabilizer and its projective action kernel give the same quotient
degree $192/8=24$.

\subsection{Extension of the canonical class}

Let \(H\) denote the pullback of a projective hyperplane to \(Y\), and
let \(E_j\) denote the exceptional divisors. Let \(\Delta\) be the
divisor with coefficient \(1/3\) along each strict transform of a
projectivized reflecting hyperplane of \(G_{29}\) and coefficient
\(2/3\) along each \(E_j\), as in
\cite[\S5.3]{DerauxVolumes}. Thus
\[
 A=K_Y+\Delta
   =\frac{28}{3}H-\frac43\sum_{j=1}^{20}E_j .
\]
Its restriction to \(S\), with plane class \(h\) and exceptional
curves \(\varepsilon_j\), is
\[
 A|_S=\frac{28}{3}h-\frac43\sum_{j=1}^{6}\varepsilon_j .
\]
For equal holomorphic sectional curvature, the canonical curvature
of a geodesic $\B^2$ is $3/4$ of the restricted canonical curvature
of $\B^3$. Put
\[
 B=\frac34 A|_S=7h-\sum_{j=1}^{6}\varepsilon_j,\qquad B^2=49-6=43.
\]
To compute the Euler characteristic from $B^2$, we extend the
curvature identity across the cusps.

\begin{lem}\label{lem:canonical-extension}
Let $\overline X=Y/PG_{29}$, and let $D^\nu$ be the normalization of
the compact mirror image in $\overline X$. Write $D^{\circ,\nu}$
for its open part over $Y^0/PG_{29}$. Let $\mathcal K$ be the
ambient orbifold canonical $\Q$-Cartier class on $\overline X$, whose
pullback to $Y$ is $A$. Choose a neat finite cover $F$ of
$\Lambda_{29}\backslash M$ induced from a neat normal ambient cover.
Its smooth toroidal compactification $\overline F$, with boundary
$C=\coprod_j C_j$, has a morphism
$f\colon\overline F\to D^\nu\to\overline X$ extending the open map.
Each $C_j$ is contracted, and
\[
 c_1(K_{\overline F}+C)=\frac34c_1(f^*\mathcal K)
 \quad\text{in }H^2(\overline F,\R).
\]
Moreover, $(K_{\overline F}+C)^2=3\chi(F)$.
\end{lem}

\begin{proof}
The toroidal and minimal compactifications of $F$ exist by
\cite[Theorem~1]{Mok}. The boundary consists of disjoint elliptic
curves with negative self-intersection, contracted to finitely many
points in the normal minimal compactification.
The finite map $F\to D^{\circ,\nu}$ is algebraic. Indeed, it is an
unramified finite cover away from the branch and singular strata,
so Riemann existence applies there; normalization extends it across
the omitted strata. Normalizing $D^\nu$ in the resulting function
field gives a normal projective compactification $T$ of $F$ with
finite boundary.

To identify $T$ with the minimal compactification, close the graph of
the identity on $F$ in the product of the two
compactifications. Any added point has both coordinates in the
finite boundaries. If one coordinate were in $F$, continuity of the
identity on $F$ would put the other there as well.
Remmert--Stein extends the graph across this finite set. Its two
projections are finite and birational, and hence are isomorphisms
by normality. Composing the toroidal contraction with the resulting
finite map to $D^\nu$ gives $f$.

On $F$, the ratio of canonical curvatures is $3/4$. The cohomology
sequence with support in the boundary therefore gives real numbers
$b_j$ such that
\[
 c_1(K_{\overline F}+C)-\frac34c_1(f^*\mathcal K)
       =\sum_j b_j[C_j].
\]
By adjunction, $(K_{\overline F}+C)\cdot C_j=0$. Also
$f^*\mathcal K\cdot C_j=0$, since $f$ contracts $C_j$.
The curves are disjoint and $C_j^2<0$, so intersection with $C_j$
gives $b_j=0$, proving the equality of Chern classes.

For the proportionality identity, apply Mumford's Chern-form theorem
\cite[Theorem~1.4]{Mumford77}. In Mok's cusp coordinates
\cite[Proposition~1 and the preceding formulas~(6)--(8)]{Mok},
write $u=-\log\|t\|$ in the cusp line bundle.
The squared norm of a logarithmic canonical frame is proportional
to $u^3$; its Chern connection and curvature have Poincar\'e growth.
Thus the metric is good in Mumford's sense, and its Chern square
integrates to $(K_{\overline F}+C)^2$. This local argument does not
require arithmeticity.
On the truncation $u=R$, the boundary term in Chern--Gauss--Bonnet
vanishes as $R\to\infty$ because
curvature and the shape operator are bounded in a horospherical
orthonormal frame, while the boundary volume is $O(R^{-2})$.
The two-ball curvature identity $c_1^2=3c_2$ now gives
$(K_{\overline F}+C)^2=3\chi(F)$.
\end{proof}

Let $m$ be the degree of $F\to\Lambda_{29}\backslash M$, and let
$\mathcal K_D$ be the restriction of $\mathcal K$ to $D^\nu$.
The projection formula gives
\[
 (f^*\mathcal K)^2=m\mathcal K_D^2,\qquad
 (A|_S)^2=24\mathcal K_D^2.
\]
Lemma~\ref{lem:canonical-extension} and multiplicativity under the
neat cover yield
\begin{equation}\label{eq:chi}
 \chiorb(\Lambda_{29})=\frac{\chi(F)}m
 =\frac{1}{3\cdot24}\left(\frac34\right)^2(A|_S)^2
 =\frac{43}{72}.
\end{equation}
Here $24$ is the degree of the coarse finite map from $S$, rather
than the degree of an orbifold covering. The order-three normal
reflection has already been divided out in $\Lambda_{29}$.

\subsection{Divisor coefficients on \texorpdfstring{$S$}{S}}

We now compute the divisor on $S$ whose log canonical class is $B$.
In the mirror coordinates above, put
\[
 u=iy_1-(1+i)y_2,\qquad v=y_1+(1+i)y_3,\qquad
 w=y_1-(1+i)y_2+(1+i)y_3.
\]
This change of coordinates has determinant $-2$. The twenty-one
restricted lines are
\begin{equation}\label{eq:restricted-lines}
\begin{aligned}
 A:\quad &iu\pm v\pm w=0,\quad u\pm iv\pm w=0,\quad
             u\pm v\pm iw=0,\\
 B:\quad &u\pm v=0,\quad u\pm w=0,\quad v\pm w=0,\\
 C:\quad &u=0,\quad v=0,\quad w=0,
\end{aligned}
\end{equation}
where the signs are independent.
The density is the pulled-back cone
angle divided by $2\pi$, and the boundary coefficient is one minus
the density. Since $S^0$ maps finitely to the mirror orbifold,
the density on $S$ need not be the reciprocal of an integer.
\begin{center}
\begin{tabular}{ccccc}
\toprule
Type & Ambient multiplicity & Number & Boundary coefficient & Density\\
\midrule
\(A\)&3&12&\(1/2\)&\(1/2\)\\
\(B\)&2&6&\(1/3\)&\(2/3\)\\
\(C\)&4&3&\(2/3\)&\(1/3\)\\
\bottomrule
\end{tabular}
\end{center}
These densities are ratios of local mirror stabilizer orders for
the finite model and the deformed orbifold, after removing the
pointwise kernels. Along a line of type
\(A\), the finite transversal group is \(S_3\), whose effective
mirror stabilizer is trivial. Its \(p=3\) counterpart has order
24 and four mirrors of order three, giving effective order \(2\).
Along a line of type \(B\), the effective stabilizer orders
are \(2\) and \(3\), respectively.
Along a line of type \(C\), the finite dihedral group has two
mirror orbits of size two,
so its effective stabilizer order is \(2\). The deformed group has order
72 and two mirror orbits of size four, hence effective order \(6\).
The local group orders are listed in Table~\ref{tab:local-groups}.

The finite local group at each blown-up point is
\(G(4,4,3)\). For its mirror \(x_1=x_2\), the effective mirror
stabilizer is \(\mu_4\), acting on tangent coordinates as
\(\operatorname{diag}(\lambda,\lambda^{-2})\).
Its action on the exceptional \(\mathbf P^1\) has multiplier
\(\lambda^3\), hence is faithful. Thus the finite quotient map is
generically unramified along each exceptional curve.
In the deformed group, put
\[
 B_{\mathrm{exc}}=(R_4R_3R_2)^8.
\]
The matrix $B_{\mathrm{exc}}$ has order three, acts as \(\omega^2 I\) on
\(\operatorname{span}(e_2,e_3,e_4)\), and is the identity on its
orthogonal complement. Thus \(\omega B_{\mathrm{exc}}\) fixes that
three-dimensional subspace pointwise and is an order-three reflection.
Its hyperplane has signature \((2,1)\) and corresponds to the
exceptional divisor. Since the ambient cone order is three and
the finite quotient is generically unramified along the exceptional
curve, its boundary coefficient is $2/3$ and its density is $1/3$.

Let $D_S$ be the divisor with these coefficients on the strict
transforms of the restricted lines and on the exceptional curves.
The total weighted line degree is
\[
 12\cdot\frac12+6\cdot\frac13+3\cdot\frac23=10.
\]
Each blown point lies on four \(A\)-lines and one \(C\)-line, with
incident weight \(8/3\). Its coefficient in the log canonical
class after blowup is
\[
 1+\frac23-\frac83=-1.
\]
Consequently
\[
 K_S+D_S=7h-\sum_{j=1}^{6}\varepsilon_j=B.
\]
At each of the nine cusp points the incident weight is \(2\).
Blowing up a cusp and assigning coefficient \(1\) to its boundary
exceptional has discrepancy \(1+1-2=0\).
Thus this ordinary blowup of the log pair preserves the log canonical
square $43$. Its exceptional curve is rational. The elliptic
boundary curves occur in the toroidal compactification of a neat
finite cover, whose canonical extension is justified by
Lemma~\ref{lem:canonical-extension}.

\section{The orbifold Euler characteristic}
\label{sec:euler}

We check \eqref{eq:chi} by pulling Satake's weighted Euler
characteristic formula \cite[\S3.3]{Satake57} back to the strata on $S$.
We write $e$ for compactly supported Euler characteristic.
The local group orders are computed from the matrices
\eqref{eq:ambient-rows} by exact enumeration in $\Q(i,\omega)$.
In Table~\ref{tab:local-groups}, the centralizer is that of the
first listed reflection. Dividing the centralizer order by three
removes the pointwise kernel and gives the effective mirror
stabilizer order in the last column. The two groups of order $648$ correspond to the
two rank-three orbits.
\begin{table}[htbp]
\centering
\begin{tabular}{lrrr}
\toprule
Generators & Order & Centralizer order & Effective mirror order\\
\midrule
$R_1,R_3$ &9&9&3\\
$R_1,R_2$ &24&6&2\\
$R_2,R_4$ &72&18&6\\
$R_1,R_2,R_3$ &648&54&18\\
$R_1,R_2,R_3^{-1}R_4R_3$ &648&54&18\\
\bottomrule
\end{tabular}
\caption{Local groups and their actions on a mirror.}
\label{tab:local-groups}
\end{table}
For each row, with generators $g_1,\ldots,g_s$, set
\[
 \mathcal G_0=\{I\},\qquad
 \mathcal G_{n+1}=\mathcal G_n\cup
       \bigcup_{j=1}^{s}g_j\mathcal G_n.
\]
All operations use rational coefficients in the basis
$1,\sqrt3,i,i\sqrt3$ of $\Q(i,\omega)$, so equality is tested exactly.
The sets stabilize at the cardinalities in the second column.
Every element is a generator word, and stability under the
generators proves that the final set is the whole group.
Counting its elements $g$ satisfying $gg_1=g_1g$ gives the third
column.

A unitary transformation preserving the mirror also preserves its
orthogonal polar line and commutes with $g_1$. Hence the setwise
mirror stabilizer equals the centralizer of $g_1$.
There are at most three generating reflections in dimension four,
so they have a nonzero common fixed vector in $\ker(g_1-I)$.
Hence a centralizer element acting scalarly on this
hyperplane fixes it pointwise. It preserves the orthogonal polar
line, on which its eigenvalue is its determinant and belongs to
$\langle\omega\rangle$. Thus the kernel of the projective mirror
action is exactly $\langle g_1\rangle$, of order three.

The line equations \eqref{eq:restricted-lines} are preserved by
coordinate permutations and sign changes. The intersection points
are obtained from the following representatives $[u:v:w]$ by these
operations. At a point on $a,b,c$ lines of types $A,B,C$, respectively,
the ambient multiplicity is $1+2a+b+3c$.
\begin{center}
\begin{tabular}{crrrrrl}
\toprule
Representative & $A$ & $B$ & $C$ & Multiplicity & Number & Treatment\\
\midrule
$[1:1:1]$&0&3&0&4&4&retain\\
$[1:1:2i]$&1&1&0&4&12&retain\\
$[1:1:1\pm i]$&2&1&0&6&24&retain\\
$[1:0:0]$&0&2&2&9&3&delete cusp\\
$[1:1:0]$&2&1&1&9&6&delete cusp\\
$[1:i:0]$&4&0&1&12&6&blow up\\
\bottomrule
\end{tabular}
\end{center}
The two signs in the third row give disjoint orbits of twelve points.
These points account for
\[
 4\cdot3+12+24\cdot3+3\cdot6+6\cdot6+6\cdot10
 =210=\binom{21}{2}
\]
line pairs, so the list is complete. The same symmetries are
transitive on each line type. Thus each line of type $A,B,C$ contains,
respectively,
\[
 \frac{12+48+12+24}{12}=8,\qquad
 \frac{12+12+24+6+6}{6}=10,\qquad
 \frac{6+6+6}{3}=6
\]
marked points. Blowing up a point replaces it on each incident line
by the intersection with the exceptional curve.
After the six blowups there are \(79\) marked points, including the
nine cusp points. The open portions of the curves of types $A,B,C$
and of the exceptional curves have Euler characteristics
$-6,-8,-4,-3$, respectively.
Let \(S^\circ\) be the complement of the strict transforms of the
restricted lines and the exceptional curves. Since \(e(S)=9\),
\[
 e(S^\circ)
 =9-\left(12(-6)+6(-8)+3(-4)+6(-3)\right)-79=80.
\]
The weighted open-curve contribution is
\begin{align*}
 &12(-6)\frac12+6(-8)\frac23
       +3(-4)\frac13+6(-3)\frac13\\
 &\hspace{40mm}=-36-32-4-6=-78.
\end{align*}

At a retained point, the density is the ratio of the finite and
deformed effective local group orders. The deleted cusps contribute
no point strata; their zero densities in the table record this deletion.
\begin{center}
\begin{tabular}{lrr}
\toprule
Point type & Number & Density\\
\midrule
Three \(B\)-lines &4&\(1/4\)\\
One \(A\)- and one \(B\)-line &12&\(1/3\)\\
Two \(A\)- and one \(B\)-line &24&\(1/9\)\\
\(A\)-line meeting an exceptional &24&\(1/6\)\\
\(C\)-line meeting an exceptional &6&\(1/9\)\\
Cusp &9&\(0\)\\
\bottomrule
\end{tabular}
\end{center}
The first row is represented by $L_{134}$, with local groups
$\langle r_1\rangle\times\langle r_3,r_4\rangle$ and
$\langle R_1\rangle\times\langle R_3,R_4\rangle$.
Removing the first factors gives effective orders $6$ and $24$,
hence density $6/24=1/4$. The second density is the product
\((1/2)(2/3)=1/3\). For the third, the finite rank-three group is
\(S_4\), with effective mirror stabilizer of order two, whereas
the deformed effective mirror stabilizer has order $18$, as in
Table~\ref{tab:local-groups}. The density is therefore
\(2/18=1/9\).

At an intersection with an exceptional curve, work in the smooth
local chart of $S$ before taking the residual finite quotient.
The two curves meet normally, so choose coordinates with the
exceptional curve $x=0$ and the other curve $y=0$.
Their complement in a small bidisk has
fundamental group $\Z^2$, generated by the two meridians.
The cone densities above give holonomy orders $(3,2)$ at an
intersection with an $A$-line and $(3,3)$ at an intersection with a $C$-line.
Taking roots of these orders therefore kills all local monodromy.
The developing map becomes single-valued and extends holomorphically
across the axes by boundedness. It takes values in the ball, since an
interior point mapping to its boundary would force the map to be
constant by the maximum principle. It is a local biholomorphism off
the axes and at their generic points by the cone-order calculation.
Its Jacobian cannot vanish only at the remaining corner, because
the zero set of a nonzero holomorphic function is a divisor. The developing
map is therefore a local biholomorphism on the whole root cover.

At the four \(A\)-directions on each exceptional curve, the finite
effective point group is trivial, so the local orbifold group has
order $3\cdot2=6$ and density $1/6$.
At the $C$-direction, the residual finite group has order four,
giving order $4\cdot3\cdot3=36$ and density $4/36=1/9$.

Besides the $C$-direction, the finite $\mu_4$-action on each
exceptional curve has one other fixed point, away from the restricted
lines. In local coordinates with exceptional curve \(x=0\), its action is
\((x,t)\mapsto(\lambda x,\lambda t)\).
The normal root \(x=z^3\) lifts the local group to order \(12\).
Its density \(4/12=1/3\) equals that of the open exceptional curve,
so we include this point in the same stratum.

The total point contribution is
\[
 4\left(\frac14\right)+12\left(\frac13\right)
 +24\left(\frac19\right)
 +6\left(4\cdot\frac16+\frac19\right)=\frac{37}{3}.
\]
Thus the total pulled-back Euler number is
\[
 80-78+\frac{37}{3}=\frac{43}{3}.
\]
Dividing by the quotient degree \(24\) gives \(43/72\), as in
\eqref{eq:chi}.

\subsection{A Deligne--Mostow example}

The canonical-class calculation of Section~\ref{sec:volume}
also applies to a diagonal mirror of
\(\mathcal C(B_4,3,4)\).
Deraux's model \cite[\S5.6.2]{DerauxVolumes} gives the ambient class
\[
 A_0=\frac{10}{3}H-\frac43\sum_{j=1}^{4}E_j.
\]
A diagonal mirror contains two of the four blown-up points, and its
finite quotient map has degree \(8\). The restriction formula gives
\[
 \frac{1}{3\cdot8}\left(\frac34\right)^2
 \left(\left(\frac{10}{3}\right)^2
       -2\left(\frac43\right)^2\right)=\frac{17}{96}.
\]
This is the known orbifold Euler characteristic of the
Deligne--Mostow quotient for \((3,3,5,6,7)/12\) with order-two
symmetry \cite[Proposition~7.10 and Appendix~A.10]{DPP}.

\section{Comparison with known nonarithmetic lattices}
\label{sec:comparison}

\begin{thm}\label{thm:comparison}
The lattice $\Lambda_{29}$ is not commensurable with any
Deligne--Mostow lattice in $\PU(2,1)$ or with any lattice in the
twenty-two nonarithmetic commensurability classes of
\cite[Theorem~1.4 and Appendix~A.11]{DPP}.
\end{thm}

Table~\ref{tab:dpp-classes} gives one representative of each of these
twenty-two classes, its adjoint trace field, and the distinction from
$\Lambda_{29}$. The notation $\Gamma(p,t)$, $\mathcal S(p,\tau)$ and
$\mathcal T(p,\mathbf T)$ follows \cite{DPP}.
Put
\[
 \mu=(3,3,5,6,7)/12,\qquad \Sigma=\Z/2\Z,
 \qquad \Gamma_{\mathrm M}=\Gamma(6,1/6),
\]
where $\Sigma$ exchanges the two equal weights.
Let $\Gamma_{\mu,\Sigma}$ be the Deligne--Mostow lattice for these
weights and symmetry group. By \cite[Proposition~7.10]{DPP}, it
contains a conjugate of $\mathcal T(4,\mathbf E_2)$ with index three.
The latter has Euler characteristic $17/32$ and represents the same
commensurability class.
The first nine rows are the Deligne--Mostow classes, also listed in
\cite[Table~2]{YuZheng}. The remaining thirteen are the additional
classes of \cite{DPP16,DPP}.

\begin{table}[htbp]
\centering
\small
\begin{tabular}{lll}
\toprule
Representative & Adjoint trace field & Distinction from $\Lambda_{29}$\\
\midrule
\multicolumn{3}{l}{Deligne--Mostow classes}\\
\midrule
$\Gamma(4,1/12)$ & $\Q(\sqrt3)$ & Compact quotient\\
$\Gamma(6,1/6)$ & $\Q(\sqrt3)$ & Proposition~\ref{prop:mostow-class}\\
$\Gamma_{\mu,\Sigma}$ & $\Q(\sqrt3)$ & Proposition~\ref{prop:gaussian-class}\\
$\Gamma(3,1/12)$ & $\Q(\sqrt2,\sqrt3)$ & Adjoint trace field\\
$\Gamma(4,3/20)$ & $\Q\bigl(\sqrt{(5+\sqrt5)/2}\bigr)$ & Adjoint trace field\\
$\Gamma(5,1/5)$ & $\Q\bigl(\sqrt{(5+\sqrt5)/2}\bigr)$ & Adjoint trace field\\
$\Gamma(9,11/18)$ & $\Q(\cos(2\pi/9))$ & Adjoint trace field\\
$\Gamma(10,4/15)$ & $\Q(\cos(2\pi/15))$ & Adjoint trace field\\
$\Gamma(3,5/42)$ & $\Q(\cos(2\pi/21))$ & Adjoint trace field\\
\midrule
\multicolumn{3}{l}{Additional sporadic and Thompson classes}\\
\midrule
$\mathcal S(4,\sigma_1)$ & $\Q(\sqrt2)$ & Adjoint trace field\\
$\mathcal S(3,\sigma_5)$ & $\Q(\sqrt5)$ & Adjoint trace field\\
$\mathcal S(3,\sigma_1)$ & $\Q(\sqrt6)$ & Adjoint trace field\\
$\mathcal S(6,\sigma_1)$ & $\Q(\sqrt6)$ & Adjoint trace field\\
$\mathcal S(4,\overline{\sigma}_4)$ & $\Q(\sqrt7)$ & Adjoint trace field\\
$\mathcal S(6,\overline{\sigma}_4)$ & $\Q(\sqrt{21})$ & Adjoint trace field\\
$\mathcal S(8,\overline{\sigma}_4)$ & $\Q(\sqrt2,\sqrt7)$ & Adjoint trace field\\
$\mathcal S(4,\sigma_5)$ & $\Q(\sqrt3,\sqrt5)$ & Adjoint trace field\\
$\mathcal T(4,\mathbf S_2)$ & $\Q(\sqrt3,\sqrt5)$ & Adjoint trace field\\
$\mathcal S(12,\overline{\sigma}_4)$ & $\Q(\sqrt3,\sqrt7)$ & Adjoint trace field\\
$\mathcal S(5,\overline{\sigma}_4)$ & $\Q\bigl(\sqrt{(5+\sqrt5)/14}\bigr)$ & Adjoint trace field\\
$\mathcal T(5,\mathbf S_2)$ & $\Q(\cos(2\pi/15))$ & Adjoint trace field\\
$\mathcal T(3,\mathbf H_2)$ & $\Q(\cos(2\pi/15))$ & Adjoint trace field\\
\bottomrule
\end{tabular}
\caption{The twenty-two nonarithmetic classes in
\cite[Appendix~A.11]{DPP} and their comparison with $\Lambda_{29}$.}
\label{tab:dpp-classes}
\end{table}

The three classes with adjoint trace field $\Q(\sqrt3)$ require
further consideration. The quotient for $\Gamma(4,1/12)$, also represented
by the weights $(4,4,5,5,6)/12$, is compact. The quotients for
$\Gamma_{\mu,\Sigma}$ and $\Gamma_{\mathrm M}$ are not compact.
They are cases 73 and 74 of \cite[Table~2]{YuZheng}, respectively;
case 74 has weights $(4,4,4,5,7)/12$.

Put $F=\Q(\sqrt3)$, so $K=F(i)$.
The matrices in Section~\ref{sec:matrices} give
$\Lambda_{29}\subset\PU(G)(F)$. Since $\Lambda_{29}$ is Zariski
dense and has adjoint trace field $F$, this is its natural $F$-form.
The three Deligne--Mostow classes have the same $K$, since their
displayed weights have least common denominator $12$; see
\cite[\S6.3]{YuZheng}.
The extension $K/F$ and the $F$-conformal class of the Hermitian form
are commensurability invariants by
\cite[Propositions~6.4(1) and~6.8]{YuZheng}.

The form $G$ has signature $(2,1)$ at both real embeddings of $F$.
For the weights $(4,4,5,5,6)/12$, $(3,3,5,6,7)/12$ and
$(4,4,4,5,7)/12$, the embedding $\zeta_{12}\mapsto\zeta_{12}^5$
gives signatures $(1,2)$, $(1,2)$ and $(2,1)$, respectively;
at the identity embedding all three have signature $(1,2)$.
Since all four forms have dimension three and their signatures agree
up to reversal at each real place,
\cite[Proposition~6.2(1)]{YuZheng} shows that all four forms are
$F$-conformal. Thus their projective unitary groups are isomorphic
over $F$. These data do not distinguish the lattice commensurability
classes. We use compactness for $\Gamma(4,1/12)$, boundary isogeny
classes for $\Gamma_{\mu,\Sigma}$, and multiplicativity of the Euler
characteristic together with an index bound for $\Gamma_{\mathrm M}$.

\subsection{The weights \texorpdfstring{$(3,3,5,6,7)/12$}{(3,3,5,6,7)/12}}

\begin{prop}\label{prop:gaussian-class}
The lattices $\Lambda_{29}$ and $\Gamma_{\mu,\Sigma}$ are not commensurable.
\end{prop}

\begin{proof}
Let $\Gamma\subset\PU(2,1)$ be a lattice such that
$\Gamma\backslash\B^2$ is not compact. Write $\mathcal A(\Gamma)$
for the set of isogeny classes of the elliptic boundary curves of a
neat finite cover. Passing to a finite cover replaces each boundary
curve by an isogenous curve, so $\mathcal A(\Gamma)$ is a
commensurability invariant \cite[Proposition~8.1]{YuZheng}.

For $\mu=(3,3,5,6,7)/12$, the unique cusp corresponds to the partition
$(5,7)\mid(3,3,6)$ of the weight numerators. Its boundary elliptic curve is isogenous
to $E_i=\C/\Z[i]$ by \cite[Example~8.2 and Corollary~8.3]{YuZheng}.
The same holds after taking the finite symmetry quotient. Hence
\[
 \mathcal A(\Gamma_{\mu,\Sigma})=\{[E_i]\},\qquad
 \mathcal A(\Lambda_{29})=\{[E_\omega]\}.
\]
Here brackets denote isogeny classes, and the second equality follows
from Section~\ref{sec:matrices}.
Isogenous elliptic curves have isomorphic rational endomorphism
algebras, whereas
\[
 \operatorname{End}(E_i)\otimes\Q=\Q(i),\qquad
 \operatorname{End}(E_\omega)\otimes\Q=\Q(\sqrt{-3}).
\]
These nonisomorphic fields imply $[E_i]\ne[E_\omega]$, excluding
the commensurability class of $\Gamma_{\mu,\Sigma}$.
\end{proof}

\subsection{Mostow's \texorpdfstring{$\Gamma(6,1/6)$}{Gamma(6,1/6)}}

\begin{prop}\label{prop:mostow-class}
The lattices $\Lambda_{29}$ and $\Gamma_{\mathrm M}$ are not commensurable.
\end{prop}

\begin{proof}
For a finite-index subgroup $\Gamma'\subset\Gamma$, the corresponding
orbifold covering gives
\[
 \chiorb(\Gamma')=[\Gamma:\Gamma']\chiorb(\Gamma).
\]
Thus $\chiorb(\Gamma)$ itself is not a commensurability invariant.
We follow the proof of \cite[Proposition~7.9]{DPP}, combining this
formula with a bound on the index in a lattice overgroup.

The lattice $\Gamma_{\mathrm M}$ has one cusp and
$\chiorb(\Gamma_{\mathrm M})=11/144$, as in the proof of
\cite[Proposition~7.9(3)]{DPP}.
Let $\widehat\Gamma\subset\PU(2,1)$ be any lattice containing
$\Gamma_{\mathrm M}$, and put $d=[\widehat\Gamma:\Gamma_{\mathrm M}]$.
This index is finite because both groups are lattices.
Fix a cusp point and let $P_{\mathrm M}$ and $\widehat P$ be its
stabilizers in $\Gamma_{\mathrm M}$ and $\widehat\Gamma$, respectively.
Since $\Gamma_{\mathrm M}$ has one cusp, the action of
$\Gamma_{\mathrm M}$ on the cusp points of $\widehat\Gamma$ is transitive.
Hence $[\widehat P:P_{\mathrm M}]=d$.

Let $B_{\mathrm M}$ be the largest open horoball precisely invariant under
$P_{\mathrm M}$ in $\Gamma_{\mathrm M}$, and write
$V_{\mathrm M}=\operatorname{Vol}(P_{\mathrm M}\backslash B_{\mathrm M})$,
using the volume normalization of \cite[\S7.2.2]{DPP}.
Proposition~7.16 of \cite{DPP}, using Parker's cusp-volume bound
\cite[Theorem~4.1]{Parker98}, gives $V_{\mathrm M}\ge d/4$.
Indeed, a cusp neighborhood for $\widehat\Gamma$ of volume at least
$1/4$ lifts to one of volume at least $d/4$ for $\Gamma_{\mathrm M}$.
The upper bound in \cite[Corollary~7.15]{DPP}, evaluated for
$\Gamma_{\mathrm M}$ in item~(3) following Proposition~7.16, gives
\[
 \frac d4\le V_{\mathrm M}\le\frac{2+\sqrt3}{4}.
\]
Consequently
\begin{equation}\label{eq:index-bound}
 [\widehat\Gamma:\Gamma_{\mathrm M}]
 \le 2+\sqrt3<4
\end{equation}
for every such lattice overgroup. Neither normality nor
torsion-freeness is required. Although Proposition~7.16 is stated in
$\mathrm{SU}(2,1)$, taking inverse images under
$\mathrm{SU}(2,1)\to\PU(2,1)$ preserves the index and the effective
cusp quotients. Thus the bound applies to the projective groups here.

Suppose that $\Lambda_{29}$ and $\Gamma_{\mathrm M}$ are
commensurable. After conjugating $\Lambda_{29}$, their commensurators
coincide. Put
\[
 \mathcal C=\Comm_{\PU(2,1)}(\Gamma_{\mathrm M}).
\]
Nonarithmeticity implies that $\mathcal C$ is a lattice containing
both groups with finite index, by Margulis' commensurator theorem
\cite[Theorem~IX.1.13]{Margulis91}; see also
\cite[Theorem~6.7 and Proposition~6.8]{DPP}.
Set
\[
 d_\Lambda=[\mathcal C:\Lambda_{29}],\qquad
 d_{\mathrm M}=[\mathcal C:\Gamma_{\mathrm M}]
\]
and apply multiplicativity to obtain
\[
 \chiorb(\mathcal C)=\frac{43}{72d_\Lambda}
                  =\frac{11}{144d_{\mathrm M}},
 \qquad
 86d_{\mathrm M}=11d_\Lambda.
\]
It follows that $11\mid d_{\mathrm M}$, in contradiction with
$1\le d_{\mathrm M}\le3$ from \eqref{eq:index-bound}.
\end{proof}

\begin{proof}[Proof of Theorem~\ref{thm:comparison}]
Section~\ref{sec:matrices} shows that $\Lambda_{29}$ is nonarithmetic
with adjoint trace field $\Q(\sqrt3)$, and that
$\Lambda_{29}\backslash\B^2$ is not compact.
Vinberg's theorem \cite[Theorems~1 and~3]{Vinberg71} makes the adjoint
trace field a commensurability invariant for Zariski-dense groups,
in particular for lattices. Thus the nineteen
rows marked ``Adjoint trace field'' in Table~\ref{tab:dpp-classes}
are excluded. Since $\Gamma(4,1/12)\backslash\B^2$ is compact,
its commensurability class is also excluded. Propositions~\ref{prop:gaussian-class}
and~\ref{prop:mostow-class} exclude the two remaining classes.

Every nonarithmetic Deligne--Mostow lattice in $\PU(2,1)$ belongs to
one of the first nine classes in the table. An arithmetic
Deligne--Mostow lattice cannot be commensurable with $\Lambda_{29}$,
since arithmeticity is a commensurability invariant. This proves
the theorem and the remaining assertion of Theorem~\ref{thm:main}.
\end{proof}

\begin{rmk}
The group $\mathcal S(6,\sigma_1)$ also has Euler characteristic
$43/72$, but its adjoint trace field is $\Q(\sqrt6)$
\cite[Appendix~A.1]{DPP}. It is therefore among the nineteen classes
excluded by the trace field.
Yu--Zheng \cite[\S8]{YuZheng} ask whether the
$(4,4,4,5,7)/12$ quotient occurs as a subball of the $G_{29}$ quotient.
Theorem~\ref{thm:main} excludes its commensurability with
$\Lambda_{29}\backslash M$. It leaves open an embedding as a different
totally geodesic surface.
\end{rmk}

\section*{Acknowledgements}
The author is grateful to AI models, especially ChatGPT for valuable assistance with the exploration of this topic and preparation of the paper. The author takes full responsibility for its content and correctness.

\bibliography{reference}

\Addresses
\end{document}